\documentclass[secnum]{jmku26}%%% eqno resets every section 
\usepackage{amsmath}
\usepackage{amssymb}
\usepackage{graphicx}

\usepackage{upref}

\newtheorem{lmm}{Lemma}[section]

\newtheorem{prp}{Proposition}[section]

\newtheorem{thm}{Theorem}[section]

\theoremstyle{definition}

\theoremstyle{remark}

\title{Classification of regular and non-degenerate projectively Anosov flows
 on three-dimensional manifolds}
\TitleHead{Regular and non-degenerate projectively Anosov flows}
\author{Masayuki ASAOKA\footnote{Department
 of Mathematics, Kyoto University, Kyoto 606-8502, Japan}}
\AuthorHead{Masayuki ASAOKA}%     %optional but recommended
\VolumeNo{00-0}%       %editors must define this.
\YearNo{0000}%         %editors must define this.
\PagesNo{000--000}%    %editors must define this.
\classificationB{37D30 (Primary); 57R30 (Secondary)}
\def\RR{\mathbb{R}}

\def\TT{\mathbb{T}}

\def\Int{\mbox{\rm{Int}}\;}

\def\Per {\mbox{\rm{Per}}}
\def\PA {\mbox{$\mathbb{P}$\rm{A}}}

\def\cO{{\mathcal O}}
\def\cF{{\mathcal F}}

\def\-{{\setminus}}
\def\ra{{\rightarrow}}
\def\st{{\;|\;}}

\def\del{\partial}
\newcommand\cl[1]{\overline{{#1}}}

\begin{document}
\maketitle

\begin{abstract}
We give a classification of
 $C^2$-regular and non-degenerate projectively Anosov flows
 on three-dimensional manifolds.
More precisely, we prove that such a flow
 must be either an Anosov flow
 or represented as a finite union of $\TT^2\times I$-models.
\end{abstract}

%
%   Introduction
%
\section{Introduction}
Mitsumatsu \cite{Mi}, and Eliashberg and Thurston \cite{ET}
 observed that any Anosov flow on a three-dimensional manifold
 induces a pair of mutually transverse
 positive and negative contact structures.
They also showed that such pairs correspond to projectively Anosov flows,
 which form a wider class than that of Anosov flows.
In \cite{ET}, Eliashberg and Thurston studied
 projectively Anosov flows, which are called conformally Anosov flows
 in their book, from the viewpoint of confoliations.

The definition of a projectively Anosov flow is as follows:
Let $\Phi=\{\Phi^t\}_{t \in \RR}$ be a flow on a
 three-dimensional manifold $M$ without stationary points.
Let $T\Phi$ denote
 the one-dimensional subbundle of the tangent bundle $TM$
 that is tangent to the flow.
The flow $\Phi$ induces a flow $\{N\Phi^t\}$ on $TM/T\Phi$.
We call a decomposition $TM=E^s + E^u$
 by two-dimensional subbundles $E^u$ and $E^s$
 {\it a $\PA$ splitting} associated with $\Phi$ if
\begin{enumerate}
 \item $E^u(z) \cap E^s(z)=T\Phi(z)$ for any $z \in M$,
 \item $D\Phi^t(E^\sigma(z))=E^\sigma(\Phi^t(z))$ for
 any $\sigma \in \{u,s\}$, $z \in M$, and $t \in \RR$, and
 \item there exist constants $C>0$ and $\lambda \in (0,1)$ such that
\begin{displaymath}
\|(N\Phi^t|_{(E^u/T\Phi)(z)})^{-1}\| \cdot 
\|N\Phi^t|_{(E^s/T\Phi)(z)}\| \leq C\lambda^t
\end{displaymath}
 for all $z \in M$ and $t>0$.
\end{enumerate} 

Note that it can be shown that each subbundle is continuous and integrable,
 and that the splitting is unique
 by the same argument for hyperbolic splittings.
A flow is called {\it a projectively Anosov flow} (or simply $\PA$ flow)
 if it admits a $\PA$ splitting.
Remark that any Anosov flow is a $\PA$ flow.
 
The subbundles $E^s$ and $E^u$ are not uniquely integrable,
 and hence, do not generate foliations in general.
We say a $\PA$ flow is ($C^r$-){\it regular}
 when both subbundles generate ($C^r$-)smooth foliations.
There are two known classes of regular $\PA$ flows.
One is the class of regular Anosov flows.
Ghys \cite{Gh} gave the complete classification of such flows.
In fact, he showed that any regular Anosov flow must be equivalent to
 either a quasi-Fuchsian flow or the suspension of
 an Anosov automorphism on the torus.
Another known class is that of flows represented by
 finite union of $\TT^2 \times I$-models given by Noda \cite{No}.
Roughly speaking,
 a $\TT^2\times I$-model is a $\PA$ flow on $\TT^2 \times [0,1]$
 which preserves the boundary tori,
 is equivalent to a linear flow on them,
 and is transverse to $\TT^2 \times \{z\}$ for any $z \in (0,1)$.
See \cite{No} for the precise definition.

A natural question is whether
 there are other regular $\PA$ flows or not.
 Noda and Tsuboi showed that no other $\PA$ flows on certain manifolds.
Their results are summarized as follows:
\begin{thm}[\cite{No},\cite{No2},\cite{NT}, and \cite{Ts}]
Let $M$ be a $\TT^2$-bundle on $S^1$ or a Seifert fibered manifold.
Then, any regular $\PA$ flow on $M$
 must be either an Anosov flow or
 represented as a finite union of $\TT^2\times I$-models.
\end{thm}

We say a dynamical system is {\it non-degenerate}
 when all periodic orbits are hyperbolic.
In this paper,
 we show that there are no new regular and non-degenerate $\PA$ flows
 on {\it any} three-dimensional manifold.
\begin{thm}
\label{thm:main}
A $C^2$-regular and non-degenerate $\PA$ flow on
 a connected and closed three-dimensional manifold
 must be either an Anosov flow or
 represented as a finite union of $\TT^2\times I$ models.
\end{thm}

The proof of Theorem \ref{thm:main} is divided into three parts.
In Subsection \ref{sec:annulus},
 we review the stability of semi-proper annular leaves
 for a $C^2$ codimension-one foliation.
Subsection \ref{sec:PA} is the main step of the proof.
We show that any regular and non-degenerate $\PA$ flow
 with a periodic orbit is Anosov.
In Subsection \ref{sec:tori},
 we show that any regular and non-degenerate $\PA$ flow
 without a periodic orbit is represented by a finite union of
 $\TT^2\times I$-models.
It is an easy consequence of the results of Arroyo and
 Rodriguez Hertz \cite{AR} and the classification by Noda.
 
\paragraph{Acknowledgment}
The author is partially supported by Grant-in-Aid
 for Scientific Research (No. 14740045),
 Japan Society for Promotion of Science, Japan.
He would like to thank Takashi Inaba
 who showed me the proof of Proposition \ref{prop:annulus},
 and Takeo Noda for many fruitful discussions.
%He is also thankful to an anonymous referee,
% who gave many suggestions to improve this paper.

%
%  Proof of Main Theorem
%
\section{Proof of Theorem \ref{thm:main}}

%
%  Semi-proper annular leaves of codimension-one foliations
%
\subsection{Semi-proper annular leaves of codimension-one foliations}
\label{sec:annulus}
In this subsection, we review the stability of semi-proper annular leaves
 for $C^2$ codimension-one foliation on three-dimensional manifolds.

For a foliation $\cF$, let $\cF(z)$ denote the leaf through a point $z$.
Recall that a leaf $L$
 of transversely orientable codimension-one foliation $\cF$ 
 is called {\it semi-proper}
 if $L$ does not accumulate to itself at lease from one side,
 and is called {\it proper}
 if $L$ does not accumulate to itself from either side.

\begin{prp}
\label{prop:annulus}
Let $M$ be a closed three-manifold
 and $\cF$ a $C^2$ codimension-one foliation on $M$.
Suppose that a leaf $L$ of $\cF$ is semi-proper
 and is homeomorphic to $S^1 \times \RR$.
Then, $L$ is a proper leaf with trivial holonomy.
\end{prp}
\begin{proof}
Suppose that $L$ is a semi-proper but not proper leaf.
By the level theory of Cantwell and Conlon,
 $L$ is at finite level,
 and hence, it is contained in an exceptional local minimal set $X$.
See Lemma 8.3.23 and Theorem 8.3.11 of \cite{CC}, for instance.
However,
 Duminy's theorem(see Theorem 1.1 of \cite{CC2}) asserts that
 the end of any semi-proper leaf of an exceptional local minimal set
 must be a Cantor set.
It contradicts that $L$ is homeomorphic to $S^1 \times \RR$.
Therefore, $L$ is a proper leaf of $\cF$.
A theorem of Cantwell and Conlon \cite[Theorem 1]{CC3}
 on the stability of ends of proper leaves
 implies that the leaf $L$ has trivial holonomy.
\end{proof}
 
%
%   regular and non-degenerate $\PA$ flows
%
\subsection{Regular and non-degenerate $\PA$ flows}
\label{sec:PA}
The main aim of this subsection is to show the following proposition.
\begin{prp}
\label{prop:Anosov}
Any $C^2$-regular and non-generate $\PA$ flow
 with a periodic orbit is Anosov. 
\end{prp}

Let $\Phi=\{\Phi^t\}$ be a $\PA$ flow
 on a closed three-dimensional manifold $M$
 and $TM=E^s +E^u$ a $\PA$ splitting associated with $\Phi$.
Suppose that $\Phi$ is $C^2$-regular and non-degenerate.
Let $\Per(\Phi)$ be the set of all periodic points of $\Phi$,
 in other words, the union of all periodic orbits.
Let $\cl{\Per(\Phi)}$ denote the closure of $\Per(\Phi)$.

Let $\cF^u$ and $\cF^s$ be the $C^2$ foliations
 generated by $E^u$ and $E^s$.
Without loss of generality, we may assume that $\cF^u$ and $\cF^s$
 are transversely orientable.
For $z \in M$,
 let $\cO(z)$ denote the orbit $\{\Phi^t(z) \st t \in \RR\}$ of $z$,
 and $\cF^s(z)$ and $\cF^u(z)$ denote the leaves of $\cF^s$
 and $\cF^u$ through $z$.
We define {\it the strong unstable set} $W^{uu}(z)$
 and {\it the unstable set} $W^u(z)$ of $z \in M$ by
\begin{displaymath}
 W^{uu}(z) =
 \{z' \in M \st \lim_{t \ra \infty} d(\Phi^{-t}(z),\Phi^{-t}(z'))=0\},
\end{displaymath}
 and $W^u(z) = \bigcup_{z' \in \cO(z)} W^{uu}(z')$.
{\it The strong stable set} $W^{ss}(z)$
 and {\it the stable set} $W^s(z)$ are defined by
 $W^{ss}(z)=W^{uu}(z;\{\Phi^{-t}\})$
 and $W^s(z)=W^u(z;\{\Phi^{-t}\})$.

The key step of the proof of Proposition \ref{prop:Anosov}
 is to show that our assumptions 
 imply $W^u(z)=\cF^u(z)$ for any $z \in \cl{\Per(\Phi)}$.
We emphasize that
 a regular $\PA$ flow may not satisfy $W^u(z)=\cF^u(z)$
 for some $z \in \cl{\Per(\Phi)}$ in general.
For example, a regular $\PA$ flow may admit
 a toral leaf $T$ of $\cF^u$ consisting of periodic orbits.
In this case, it is easy to see that
 $W^u(z)$ coincides with $\cO(z)$,
 and hence, is a proper subset of $T=\cF^u(z)$ for any $z \in T$.

First, we investigate the topology of the unstable sets.
We say a $\Phi$-invariant embedded torus $T$
 is {\it irrational} if the restriction of the flow on $T$ is
 topologically conjugate to an irrational linear flow on the torus.
The following is an immediate corollary of Theorem B of \cite{AR}.
\begin{prp}
\label{prop:decomposition}
Let $\Omega(\Phi)$ be the non-wandering set of $\Phi$.
Then, there exists a mutually disjoint decomposition
 $\Omega(\Phi)=\Omega_0 \sqcup \Omega_1 \sqcup \Omega_2$ such that
\begin{enumerate}
\item $\Omega_1$ is a compact hyperbolic set of saddle type,
\item $\Omega_0$ is the union of finitely many attracting periodic orbits
 and irrational toral leaves of $\cF^u$.
\item $\Omega_2$ is the union of finitely many repelling periodic orbits
 and irrational toral leaves of $\cF^s$.
\end{enumerate}
\end{prp}
Notice that an irrational toral leaf $T$ of $\cF^u$
 is normally attracting, and hence,
 $W^u(z)=\cO(z)$ for any $z \in T$
 and $\bigcup_{z \in T}W^s(z)$ is a neighborhood of $T$
 which is homeomorphic to $\TT^2 \times \RR$.
Using the theory of hyperbolic invariant sets
 (see Chapter 9 of \cite{Sh} for example), we have
 $M =\{\bigcup_{z \in \Omega_1 \cup \Omega_2}W^u(z)\} \cup \Omega_0$
 and $M = \{\bigcup_{z \in \Omega_0 \cup \Omega_1}W^s(z)\} \cup \Omega_2$.

\begin{lmm}
\label{lemma:loop}
Let $\sigma$ be either $u$ or $s$, $p$ a repelling periodic point,
 and $V$ an annular neighborhood of $\cO(p)$ in $\cF^\sigma(p)$.
If a leaf $L$ of $\cF^s$ contains
 a simple closed curve $\gamma \subset W^u(p)$
 which is not null-homotopic in $L$,
 then $\Phi^{-t}(\gamma) \subset V$ for some $t>0$.
In particular, we have $L=\cF^s(p)$.
\end{lmm}
\begin{proof}
Since $p$ is repelling,
 there exists a neighborhood $U$ of $\cO(p)$
 such that $U$ is homeomorphic to $S^1 \times [-1,1]^2$,
 the restriction $\cF^\sigma|_U$ of $\cF^\sigma$ on $U$
 has a unique annular leaf $L_0$ with $\cO(p) \subset L_0 \subset V$,
 and other leaves of $\cF^\sigma|_U$ are homeomorphic to $\RR \times [-1,1]$.
See Figure \ref{fig:loop}.
\begin{figure}[ht]
\begin{center}
\includegraphics[scale=0.8]{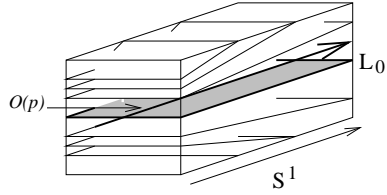}
\end{center}
\caption{The restriction of $\cF^\sigma$ on $U$}
\label{fig:loop}
\end{figure}

Suppose that $L \cap W^u(p)$ contains a simple closed curve $\gamma$
 which is not null-homotopic in $L$.
Choose a sufficiently large $t_* >0$
 so that $\Phi^{-t_*}(\gamma)$ is contained in $U$.
Since $\Phi^{-t_*}(\gamma)$ is not null-homotopic in $L$,
 we obtain $\Phi^{-t_*}(\gamma) \subset L_0 \subset V$.
\end{proof}

\begin{lmm}
\label{lemma:connected}
 $W^u(z) \cap \cF^u(z)$ is an open subset of $\cF^u(z)$
 for any $z \in M \-\Omega_0$.
\end{lmm}
\begin{proof}
The lemma is trivial if $z$ is contained in the unstable set
 of a repelling periodic point.
Hence, it is sufficient to show that
 $W^u(z)$ is tangent to $E^u$
 for any $z \in \Omega_1 \cup \Omega_2$
 which is not a repelling periodic point.
By the local unstable manifold theorem
 or the local strong unstable manifold theorem,
 there exists an injectively immersed open two-dimensional manifold $V$
 such that $\bigcap_{t>0} \Phi^{-t}(V)=\cO(z)$,
 $V \subset \Phi^t(V)$ for any $t>0$,
 and $T_{z'} V=E^u(z')$ for any $z' \in \cO(z)$.
Then, the domination property and the invariance
 of the splitting $TM=E^s+E^u$
 implies that $V$, and hence, $W^u(z)$ must be tangent to $E^u$.
\end{proof}
For $z \in M \- \Omega_0$,
 let $V^u(z)$ denote the connected component of $W^u(z) \cap \cF^u(z)$
 that contains $z$.
By the above lemma, $V^u(z)$ is an open subset of $\cF^u(z)$.

For a topological space $X$ and its subspace $Y$,
 we say a point $x \in X \- Y$ is {\it accessible} from $Y$
 if there exists a continuous map $l:[0,1] \ra X$ such that
 $l(1)=x$ and $l([0,1)) \subset Y$.
\begin{lmm}
\label{lemma:accessible}
For any $z \in M \- \Omega_0$,
 if $q \in \cF^u(z) \- V^u(z)$ is accessible from $V^u(z)$
 then it is an attracting periodic point.
\end{lmm}
\begin{proof}
By Proposition \ref{prop:decomposition}, 
 there are three possibilities:
 $q \in W^u(\Omega_1 \cup \Omega_2)$,
 $q$ is contained in an irrational toral leaf $T$ of $\cF^u$,
 or $q$ is an attracting periodic point.
If the first occurs,
 then $V^u(q) \cap V^u(z) \neq \emptyset$
 since $V^u(q)$ is an open subset of $\cF^u(z)$.
It implies that $W^u(z)= W^u(q)$.
However, it contradicts $q \not\in W^u(z)$.
The second implies $T=\cF^u(q)=\cF^u(z)$ is a subset of $\Omega_0$.
It contradicts that the assumption $z \in M \- \Omega_0$.
\end{proof}

\begin{lmm}
\label{lemma:repelling}
For a repelling periodic point $p$ of $\Phi$,
\begin{enumerate}
 \item $V^u(z)=\cF^u(z)$ for any $z \in W^u(p) \- V^u(p)$, and
 \item if $V^u(p) \neq \cF^u(p)$, then there exists
 an attracting periodic point $q \in \cF^u(p)$
 and an embedded closed annulus $A \subset \cF^u(p)$
 such that $\del A=\cO(p) \cup \cO(q)$ and
 $\Int A \subset V^u(p) \cap W^s(q)$.
\end{enumerate}
\end{lmm}
\begin{proof}
Suppose $V^u(z) \neq \cF^u(z)$ for $z \in W^u(p)$.
Take a point $q \in \cF^u(z) \- V^u(z)$ which is accessible from $V^u(z)$.
By Lemma \ref{lemma:accessible}, $q$ is an attracting periodic point.
Hence, we can take a compact annulus $A_0 \subset \cF^u(z)=\cF^u(q)$
 so that  $A_0 \- \cO(q) \subset V^u(z)$,
 where $\cO(q)$ is a boundary component of $A_0$
 and the other boundary component $\gamma$ is transverse to the flow.
By the Poincar\'e-Bendixon theorem,
 $\gamma$ is not null-homotopic in $\cF^u(z)$.
Take an annular neighborhood $V$ of $\cO(p)$ in $V^u(p)$
 and apply Lemma \ref{lemma:loop} to $\sigma=u$, $L=\cF^u(z)$, and $V$.
Then, we obtain $t>0$ satisfying $\Phi^{-t}(\gamma) \subset V$.
In particular, we have $\gamma \subset V^u(p)$,
 and hence, $V^u(z)=V^u(p)$.
It is easy to find the required closed annulus
 in $\Phi^{-t}(A_0) \cup V$.
\end{proof}

The next proposition is the key step of the proof.
\begin{prp}
\label{prop:repelling}
The flow $\Phi$ has neither repelling nor attracting periodic points.
In particular,
 $\cF^u(z)=W^u(z)$ and $\cF^s(z)=W^s(z)$ for any $z \in \Omega_1$.
\end{prp}
\begin{proof}
We show that $\Phi$ has no repelling periodic points.
It also implies the non-existence of attracting periodic points
 once we replace $\Phi=\{\Phi^t\}$ by $\{\Phi^{-t}\}$.
Then, the latter assertion follows from Lemma \ref{lemma:accessible}
 and the fact that $W^u(z)$ is an immersed manifold tangent to $E^u$
 for any $z \in \Omega_1$.

Assume that there exists a repelling periodic point $p$.
Take an annular neighborhood $V_0$ of $p$ in $\cF^u(p)$ such that
 $\Phi^{-t}(V_0) \subset V_0$ for any $t>0$
 and $\bigcap_{t>0} \Phi^{-t}(V_0)=\cO(p)$.
Note that $V^u(p)=\bigcup_{t>0} \Phi^t(V_0)$.
We also take a neighborhood $U$ of $\cO(p)$ in $M$
 such that $U \cap \del V_0 =\emptyset$
 and $\Phi^{-t}(U) \subset U$ for any $t >0$.
We claim $U \cap V^u(p) =U \cap V_0$.
In fact, suppose $\Phi^{t_0}(z) \in U$ for $z \in V_0$
 and $t_0 >0$.
Then, we have $\Phi^t(z) \in U$ for any $t \in [0,t_0]$.
If $\Phi^{t_0}(z) \not\in V_0$,
 then there exists $t_1 \in (0,t_0)$ such that
 $\Phi^{t_1}(z) \in \del V_0$
 since $\Phi^t(V_0) \supset V_0$ for any $t>0$.
It contradicts that $U \cap \del V_0= \emptyset$.
Hence, $\Phi^{t_0}(z)$ is a point of $V_0$.

If $\cF^u(p)=V^u(p)$,
 the leaf $\cF^u(p)$ is proper since $U \cap V^u(p)=U \cap V_0$.
However, it contradicts Proposition \ref{prop:annulus}
 since $\cF^u(p)=V^u(p)=\bigcup_{t>0}\Phi^t(V_0)$ is homeomorphic to
 $S^1 \times \RR$ but the linear holonomy along $\cO(p)$ is not trivial.

Second, we assume $\cF^u(p) \neq V^u(p)$.
By Lemma \ref{lemma:repelling},
 there exist an attracting periodic point $q$
 and an embedded compact annulus $A$ such that
 $\del A= \cO(p) \cup \cO(q)$ and $\Int A \subset V^u(p) \cap W^s(q)$.
It is important to remark that the orientations of the orbits
 of $p$ and $q$ are opposite in $A$.
It is because the holonomy of $\cF^u$ along $\cO(p)$ is expanding
 and that along $\cO(q)$ is contracting.

Since $q$ is attracting,
 there exists an annular neighborhood $A^s$ of $\cO(q)$ in $\cF^s(q)$ 
 such that $\Phi^t(A^s) \subset A^s$ for any $t>0$,
 $\bigcap_{t \geq 0}\Phi^t(A^s)=\cO(q)$,
 and $\cF^u(z) \cap (U \- V_0) \neq \emptyset$
 for any $z \in A^s \- \cO(q)$.
Since $U \cap V^u(p)=U \cap V_0$,
 Lemma \ref{lemma:repelling} implies that
 $\cF^u(z') \subset W^u(p)$ for any $z' \in U \- V_0$.
Hence, we have $A^s \- \cO(q) \subset W^u(p)$.
See Figure \ref{fig:repelling}.
\begin{figure}[ht]
\begin{center}
\includegraphics[scale=1.0]{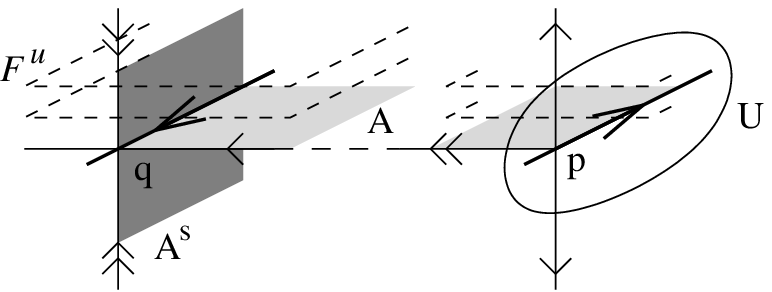}
\end{center}
\caption{Proof of Lemma \ref{lemma:repelling}.}
\label{fig:repelling}
\end{figure}
The Poincar\'e-Bendixon theorem implies that
 each boundary component of $A^s$ is not null-homotopic in $\cF^s(q)$.
Applying Lemma \ref{lemma:loop} to $\sigma=s$ and $L=\cF^s(q)$
 we obtain $\cF^s(p)=\cF^s(q)$.

Put $\lambda^\sigma_p=\|N\Phi^{T_p}|_{(E^\sigma/T\Phi)}(p)\|$
 and $\lambda^\sigma_q=\|N\Phi^{T_q}|_{(E^\sigma/T\Phi)}(q)\|$
 for $\sigma=u,s$, where $T_p$ and $T_q$ are the period of $p$
 and $q$ respectively.
Since $\cF^u(p)=\cF^u(q)$, $\cF^s(p)=\cF^s(q)$,
 and the orientations of the orbits of $p$ and $q$ are opposite,
 we have $\lambda^u_p \lambda^u_q=\lambda^s_p\lambda^s_q=1$.
On the other hand, we also have
 $\lambda^s_p<\lambda^u_p$ and $\lambda^s_q<\lambda^u_q$
 since $\Phi$ is a $\PA$ flow.
It is a contradiction.
\end{proof}

\begin{lmm}
\label{lemma:boundary}
The set $\Omega_1$ coincides with $\cl{\Per(\Phi)}$.
For $z \in \Omega_1$,
 the leaf $\cF^s(z)$ contains a periodic orbit
 if it is semi-proper.
\end{lmm}
\begin{proof}
Since $\Omega_0$ consists of normally attracting invariant tori
 and $\Omega_2$ consists of normally repelling invariant tori,
 we can take a neighborhood $U$ of $\Omega_1$ so that
 $\Omega_1=\bigcap_{t \in \RR}\Phi^t(U)$, in other words,
 $\Omega_1$ is a locally maximal hyperbolic set.
By the theory of hyperbolic invariant sets,
 we obtain $\cl{\Per(\Phi)}=\Omega_1$.
See \cite{Sh}.

The second assertion follows from Proposition \ref{prop:repelling}
 and a variant of Proposition 1 of \cite{NP} for three-dimensional flows.
\end{proof}

\begin{proof}[Proof of Proposition \ref{prop:Anosov}]
We show $M=\Omega_1$.
It implies that $\Phi$ is an Anosov flow
 since $\Omega_1$ is a hyperbolic set
 by Proposition \ref{prop:decomposition}.

By the Birkhoff-Smale theorem (see {\it e.g.} Corollary 6.5.6 of \cite{KH}),
 we have $W^u(p) \cap W^s(p) \subset \cl{\Per(\Phi)}$
 for any $p \in \Per(\Phi)$.
Proposition \ref{prop:repelling} and Lemma \ref{lemma:boundary} imply
 $\cl{\Per(\Phi)}=\Omega_1$ and
 $\cF^\sigma(z)=W^\sigma(z)$ for any $z \in \Omega_1$
 and any $\sigma =u,s$.
Since $\cF^u$ and $\cF^s$ are mutually transverse,
 we obtain that $\cF^u(p) \cap \cF^s(p) \subset \cl{\Per(\Phi)}=\Omega_1$
 for any $z \in \Omega_1$.

Suppose $\cF^u(z) \not\subset \Omega_1$
 for some $z \in \Omega_1$.
Then, there exists $z_0 \in \cF^u(z) \cap \Omega_1$
 which is accessible from $\cF^u(z) \- \Omega_1$.
It implies that $\cF^s(z_0)$ is a semi-proper leaf.
By Lemma \ref{lemma:boundary},
 $\cF^s(z_0)=W^s(z_0)$ contains a periodic point.
In particular, it is diffeomorphic to $\RR \times S^1$.
However, it contradicts Proposition \ref{prop:annulus}
 since $\cF^u(z_0)$ has the non-trivial linear holonomy
 along the periodic orbit.
Therefore, we obtain $\cF^u(z) \subset \Omega_1$
 for any $z \in \Omega_1$.
Similarly,
 we also have $\cF^s(z) \subset \Omega_1$
 for any $z \in \Omega_1$.
Since $\cF^u$ and $\cF^s$ are mutually transverse,
 it implies that $\Omega_1$ is an open subset of $M$.
Therefore, $\Omega_1=\cl{\Per(\Phi)}$ coincides with $M$
 if $\Phi$ has a periodic point.
\end{proof} 

%
%   Regular $\PA$ flows without periodic points
%
\subsection{Regular $\PA$ flows without periodic points}
\label{sec:tori}
To complete the proof of Theroem \ref{thm:main},
 we show that any regular $\PA$ flow $\Phi$ 
 without periodic points
 is represented by a finite union of $\TT^2\times I$-models.
In \cite{No},
 Noda showed that any $C^2$ regular $\PA$ flow
 on $\TT^2$-bundle over the circle
 is either an Anosov flow or is represented by a finite union of
 $\TT^2\times I$-models.
Hence, it is sufficient to show that $M$ is such a manifold.

By Proposition \ref{prop:decomposition} and Lemma \ref{lemma:boundary},
 the non-wandering set of $\Phi$ is the union of
 irrational toral leaves of $\cF^u$ and $\cF^s$.
Put $W^u(T)=\bigcup_{z \in T} W^u(z)$
 and $W^s(T)=\bigcup_{z \in T} W^s(z)$
 for such a toral leaf $T$.
Recall that if $T$ is an irrational toral leaf of $\cF^u$,
 then $W^s(T)$ is a neighborhood of $T$
 which is homeomorphic to $\TT^2 \times (-1,1)$.
Since $M=\bigcup_{z \in \Omega_0} W^s(z) \cup \Omega_2$,
 each boundary component of $W^s(T)$ is a toral leaf of $\cF^s$.
The similar holds for normally repelling tori.
Therefore, there exists a covering map
 $h:\TT^2 \times \RR \ra M$
 such that
\begin{enumerate}
\item $T_{2k}=h(\TT^2\times \{2k\})$ is an irrational toral leaf of $\cF^s$
 with $W^u(T_{2k})=h(\TT^2\times (2k-1,2k+1))$ and
\item $T_{2k+1}=h(\TT^2\times \{2k+1\})$ is
 an irrational toral leaf of $\cF^u$
 with $W^s(T_{2k+1})=h(\TT^2\times (2k,2k+2))$
\end{enumerate}
 for any integer $k$.
It implies that $M$ is a $\TT^2$-bundle over the circle.

%
%   References
%

\end{document}